\documentclass[12pt,reqno]{article}

\usepackage{amssymb}
\usepackage{amsmath}
\usepackage{amsthm}
\usepackage{amsfonts}
\usepackage{url}
\usepackage{hyperref}
\usepackage{breakurl}
\usepackage{graphicx}
\usepackage{comment}

\renewcommand{\ge}{\geqslant}
\renewcommand{\le}{\leqslant}

\newcommand{\sinc}{{\rm sinc}}

\newcommand{\R}{{\mathbb R}}

\newcommand{\raisedot}{\raisebox{0.17em}{$.$}}

\theoremstyle{plain}

\theoremstyle{definition}

\theoremstyle{remark}

\begin{document}
\bibliographystyle{plain}
\title{Jonathan Michael Borwein 1951 -- 2016:\\
Life and Legacy\\[5pt]
}
\author{Richard P. Brent\footnote{Mathematical Sciences Institute,
Australian National University, \hbox{Canberra, ACT.\ \ \ }
\hbox{Email: {\tt <Borwein@rpbrent.com>.}}}
}

\date{}

\maketitle

\begin{abstract}
Jonathan M.\ Borwein (1951--2016) was a prolific mathematician 
whose career spanned several countries
(UK, Canada, USA, Australia) and whose many interests included
analysis, optimisation, number theory, special functions, 
experimental mathematics, mathematical \hbox{finance}, mathematical education,
and visualisation. We describe his life and legacy, and give an
annotated bibliography of some of his most significant books and papers.
\end{abstract}

\section{Life and Family}				\label{sec:intro}

Jonathan (Jon) Michael Borwein was born in St Andrews, Scotland, on 20 May 1951.
He was the first of three children of 
David Borwein (1924--2021)
and Bessie Borwein
(n\'ee Flax). It was an itinerant academic family.
Both Jon's father David and his younger brother
Peter Borwein (1953--2020) are well-known mathematicians and occasional
co-authors of Jon. His mother Bessie is a former professor of anatomy.
The Borweins have an Ashkenazy Jewish background.
Jon's father was born in Lithuania, moved in 1930 with his family to
South Africa (where he met his future wife Bessie), and moved with
Bessie to the UK in 1948. There he obtained a PhD (London) 
and then a Lectureship
in St Andrews, Scotland, where Jon was born and went to school at
Madras College. The family, including Jon and his two siblings, 
moved to Ontario, Canada, in 1963.
In 1971 Jon graduated with a BA (Hons Math) from the University of
Western Ontario. It was in Ontario that Jon met his future
wife and lifelong partner Judith (n\'ee Roots).

When Jon won an Ontario Rhodes scholarship 
he returned to the UK,
where he obtained an MSc (1972) and
DPhil (Jesus College, Oxford, 1974).
He then moved back to North America, where he held various positions
at Dalhousie, CMU (Assistant/Associate Professor, 1980--82), 
Waterloo, and Simon Fraser.
These included Professor at Dalhousie (1984--1991), 
Waterloo (1991--1993), and Simon Fraser (1993--2004),
where he was founding Director of the Centre for Experimental
and Constructive Mathematics (CECM).
At the end of his time in Canada, he held
a Canada Research Chair in Collaborative Technology
at Dalhousie (2004--2009).
At different times both Jon and his father David were
presidents of the Canadian Mathematical Society.

In 2009, in the final stage of his career, Jon and Judith moved, with their
daughters Naomi and Tova and grandson Jakob, to Newcastle, Australia.  At
the University of Newcastle he became Laureate Professor of Mathematics and
founding Director of the Priority Research Centre CARMA (Computer Assisted
Research in Mathematics and its Applications).
Jon was elected a Fellow of the Australian Academy of Science in 2010. 

\begin{figure}[ht]
\begin{center}
\includegraphics[width=15em]{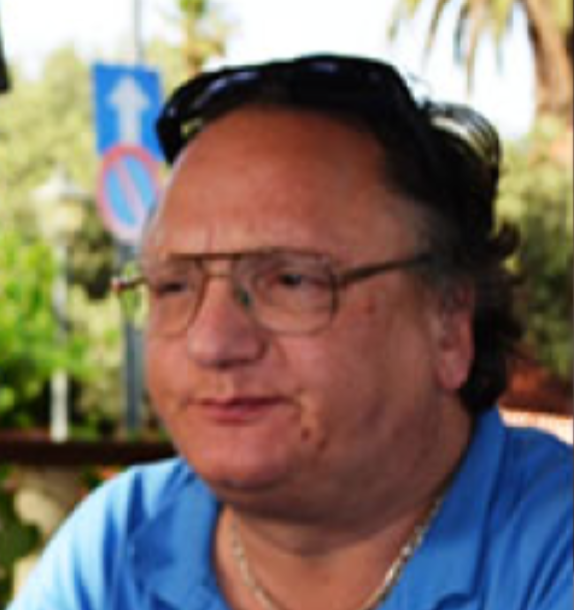}
\caption{Jon in Newcastle, c.\ 2010.}
\end{center}
\end{figure}
\vspace*{-20pt}

Jon had an outgoing personality, and loved to share his thoughts on
mathematical (and other) topics with students and colleagues.  This resulted
in him having a large number of collaborators, ranging from his father and
brother to his last PhD student.\footnote{An incomplete 
list of Jon's co-authors may be found via Google Scholar citations.}
For more on Jon's personality, see~\cite{obituary-AustMS} 
and \cite[Preface]{JBCC-proceedings}.

Jon died unexpectedly on 2 August 2016,
while on leave from Newcastle and on
a 4-month visit to Canada as Distinguished Scholar in Residence at Western
University, London, Ontario.
He was survived by his wife Judith, 
their three daughters Naomi, Rachel and Tova,
and five grandchildren, as well as by his sister Sarah, brother
Peter (d.~2020), and parents David (d.~2021) and Bessie.

\section{Mathematical Legacy}

Jon was a true ``polymath''.  While it is no longer possible for one person
to be an expert in all the diverse areas of mathematics, Jon was without
doubt an expert in several distinct areas, and very well-informed about many
other areas of mathematics and its applications.

The areas in which Jon was an expert and made notable contributions include
applied analysis, optimisation, experimental mathematics, visualisation, 
special functions, number theory, mathematical finance,
and mathematical education.  
He was a great populariser of mathematics, writing several books for a general
mathematical audience, and giving many popular talks. The slides for
some of these talks remain on his homepage~\cite{JB-homepage},
although, given the ephemeral nature of web sites, we do not know for how
much longer they will be accessible.

Jon started his mathematical career with a $1974$ DPhil in Optimisation from
Oxford, followed by a post-doctoral position in a strong functional
analysis group in Dalhousie. His early work included papers on operations
research, semi-definite programming, and integer programming.
After moving to Carnegie Mellon University in 1980, he added
variational analysis. In the eighties, by using his knowledge of Banach spaces,
he proved some results~\cite{B3} in optimisation theory 
that were surprising as they contradicted the expert opinion of the time.
Some of Jon's work in convex analysis and optimisation is described in
his book~\cite{Convex-optimization-book}, written 
with his former postdoc Adrian Lewis.

Since the mid-1980s Jon was recognised as a leader in the field of
nonlinear analysis. In later decades, while he added interests, he 
maintained activity in his earlier fields of specialisation.

Around the mid-1980s, if not earlier, Jon 
became ``sold'' on computation, aided by software such as Maple, 
as a tool for discovery in mathematics.
In his own research he was engaged in numerically
intensive optimisation, primarily the use of maximum entropy methods
(MEM). An early outcome was his highly cited paper~\cite{B4} on the
\emph{Barzilai-Borwein} optimisation algorithm.
This motivated heuristics for minimizing functions while avoiding the heavy
cost of line-search. Key to MEM is efficient numerical integration, a topic
central to Jon's research in experimental mathematics, see
for example~\cite[\S4]{closed-forms}.

From 1993 to 2004, Jon collaborated with Vancouver General Hospital's
Medical Imaging Group~\cite{Bauschke-Celler-Noll-Borwein}. 
Jon's second most-cited paper \cite{B6}\footnote{Citation data 
as at July 2021.} was relevant to this interaction.
It provided a full treatment of the algebraic iterative methods for
convex problems that are critical for turning measurements into images.
 
In 2004--2007, Jon, with collaborators David Bailey and others, completed
his first two books \cite{maths-by-experiment,maths-in-action} 
on Experimental Mathematics. These books helped to popularise experimental
mathematics, and were based largely on Jon's research papers. 
Jon continued to publish significant papers in applied experimental
mathematics, often with collaborators David Bailey, David Broadhurst, and
Richard Crandall. Also, when he moved to Newcastle in 2009, he continued to 
emphasise the importance of experimental mathematics, and founded
the research centre CARMA. Up to the time of Jon's death,
CARMA was successful in fostering the development of 
experimental mathematics within Australia.\footnote{As at 2021,
CARMA continues without Jon, but with reduced personnel and funding.}

\subsection{Nonlinear Analysis and Optimisation}	\label{subsec:analysis}

Jon's early (and some of his most significant) research was in the theory of
optimisation and related convex analysis.  A notable paper, already
mentioned above, is~\cite{B4}, which describes what is now known as the
\emph{Barzilai-Borwein algorithm} for large-scale unconstrained
optimisation. 

Suppose that we want to approximate a stationary point of a 
multi\-variate
function $F(x)$ that is differentiable in the neighbourhood of a starting
point $x_0\in \R^n$. Several optimisation methods are based on the
iteration
\[
x_{k+1} = x_k - \gamma_k\nabla F(x_k), \; k \ge 0,
\]
but differ in their choice of ``step sizes'' $\gamma_k$.
The choice used in the Barzilai-Borwein algorithm is\footnote{The paper
\cite{B4} gives two different ways of choosing $\gamma_k$;
we describe one of them here.}
\[ 
\gamma_k = \frac{(x_k-x_{k-1})^T (\nabla F(x_k)-\nabla F(x_{k-1}))}
	   {||\nabla F(x_k) - \nabla F(x_{k-1})||_2^2}\,\raisedot
\]
The motivation for this choice is that it provides a two-point approximation
to the secant equation underlying quasi-Newton
methods~\cite[\S2]{B4}. This generally gives much faster convergence
than the classical method of steepest descent, while having
comparable cost per iteration and storage requirements.

Jon's best-known contribution to optimisation theory and nonlinear analysis
is perhaps his
paper~\cite{B3} with David Preiss, which introduces a ``smooth'' variational
principle that has very broad applicability in problems of analysis.
This extends the well-known variational principle of 
Ekeland~\cite{Ekeland}.

\subsection{Experimental Mathematics}			\label{subsec:expt}

Jon was a great advocate of \emph{experimental mathematics}.  He used
advanced computer hardware and software both to discover new mathematics and
to suggest, prove, or disprove, various interesting conjectures.  We give
two examples here.  Many other examples may be found in Jon's papers and
books~\cite{maths-by-experiment,maths-in-action,crucible}.

In the first example,
an experimental computation suggested a conjecture
(subsequently proved) involving Euler's constant.  
An $n$-fold integral that arises in quantum field theory is
\[
C_n := \frac{4}{n!}\int_0^\infty\cdots\int_0^\infty
 \left(\sum_{j=1}^n (u_j + u_j^{-1})\right)^{\!-2}
 \frac{{\rm d} u_1}{u_1}\cdots\frac{{\rm d} u_n}{u_n}\,\raisedot
\]
In~\cite{Ising-class}, Jon considered the integrals $C_n$, along with
two similar (but slightly more complicated) integrals $D_n$ and $E_n$
that arise in mathematical physics, where the $D_n$ are called
\emph{Ising susceptibility integrals}. He obtained upper and lower bounds
on the $C_n, D_n$, and $E_n$ integrals for large $n$, as well as closed-form
expressions for $n \le 4$. Some examples of the latter
are $C_4 = 7\zeta(3)/12$
and $D_4 = 4\pi^2/9 - 7\zeta(3)/2 - 1/6$,
where $\zeta$ denotes the Riemann zeta-function.

Jon showed that the $C_n$ integrals can be reduced to one-dimensional
integrals involving the modified Bessel function $K_0(t)$:
\begin{equation}
								\label{eq:Cn1D}
C_n = \frac{2^n}{n!}\int_0^\infty t\,K_0^n(t)\,{\rm d}t\,.
\end{equation}
No similar reduction is known for the $D_n$ or $E_n$ integrals.

Using the representation~\eqref{eq:Cn1D} and 
tanh-sinh quadrature~\cite{tanh-sinh} to evaluate the integrals,
it was evident that the sequence $(C_n)$ is monotonic decreasing,
with a limit, say $C_\infty$, given numerically by
\[	
C_\infty = 0.63047350337438679612204019271087890435458707871273\ldots
\]
Using the Inverse Symbolic Calculator~\cite{ISC},
Jon identified $C_\infty$ with a simple closed form
\begin{equation}
C_\infty = 2e^{-2\gamma}\,,				\label{eq:Climit}
\end{equation}
where $\gamma$ is Euler's constant.
Finally, in~\cite[Theorems 1--2]{Ising-class}, Jon proved both the 
monotonicity of the sequence $(C_n)$ and convergence to the
constant~\eqref{eq:Climit}.

Some pure mathematicians might now discard all evidence of the numerical
computations, and just publish the final monotonicity and convergence
theorems. However, Jon recognised that such an approach would be 
uninformative and little short of fraudulent,
since the theorems would never have been discovered without
the preliminary numerical computations. 

Our second example is taken from a paper~\cite{surprising-sinc}
that Jon wrote with David Borwein (his father) and Robert Baillie.
The paper gives several surprising identities (and approximate identities)
involving the $\sinc$ function.
We recall that $\sinc(x)$ is defined to be $\sin(x)/x$ if $x \ne 0$,
and $1$ if $x = 0$.

In~\cite[example 1(b)]{surprising-sinc},
Jon asked the following question: for which positive integers $N$ does
\begin{equation}	
						\label{eq:sincsum1}
\frac{1}{2} + \sum_{n=1}^\infty\prod_{k=0}^N\sinc\left(\frac{n}{2k+1}\right)
  = \int_0^\infty\prod_{k=0}^N\sinc\left(\frac{x}{2k+1}\right)\,{\rm d}x\; ?
\end{equation}
Experimentation with Maple or Mathematica suggests
that~\eqref{eq:sincsum1} is an
identity for (at least) $1 \le N \le 6$. However, in a warning that it is not
always safe to extrapolate such results, Jon 
showed that \eqref{eq:sincsum1} holds for $1 \le N \le 40248$, but fails for
all $N \ge 40249$, although the difference between the left-hand and 
right-hand sides of \eqref{eq:sincsum1}
is tiny. 
Several similar, though less extreme, examples may be found
in Jon's papers~\cite{remarkable-sinc,surprising-sinc}.
 
\subsection{Number Theory and Special Functions}  \label{subsec:number_theory}

Jon was fascinated by the constant $\pi = 3.14159265358979\cdots$, 
and gave many fascinating talks on
this topic, some associated with an annual celebration of
``pi day'' on March 14th (US-style date).\footnote{At the time of writing, 
overheads for most of these talks are available at 
\url{https://www.carma.newcastle.edu.au/resources/jon/index-talks.shtml}.
}
Indeed, a talk on $\pi$ is a good way to introduce some interesting
mathematics to a general audience with a mathematical background,
or to undergraduate students. Mathematical topics that can be motivated
by $\pi$ include the concepts of irrational and transcendental numbers,
rates of convergence of series and algorithms, normality of the digits of
$\pi$ in decimal or binary (still an open problem), and Euler's famous
relation $e^{i\pi} = -1$ that connects the constants $e, \pi, i,$ and $-1$.

Much of Jon's work on $\pi$ may be found in the $1987$ book
\emph{Pi and the AGM} \cite{B12} that Jon wrote with his brother Peter.
Here \emph{AGM} refers to the arithmetic-geometric mean of Gauss
and Legendre. 

In 1975--76, the present author and (independently) Salamin published
the first known \emph{quadratically convergent} algorithm for
computing\footnote{Here ``computing $\pi$'' means approximating
$\pi$ to any desired accuracy. 
For example, Archimedes found that 
$3\frac{10}{71} < \pi < 3\frac{1}{7}$\,.
It is impossible to compute the decimal or binary representation
of $\pi$ \emph{exactly}, since $\pi$ is 
transcendental~\cite[\S11.2]{B12},
so it has an infinite, nonperiodic, decimal (or binary) representation.}
$\pi$, often called the \emph{Gauss-Legendre algorithm},
see for example~\cite{rpb252}.
In their book \emph{Pi and the AGM}, Jon and Peter Borwein gave further
quadratically (and some more rapidly) convergent algorithms. They also
discussed the complexity of computing algebraic functions and elementary
functions, the transcendence of $\pi$, the history of calculations of $\pi$,
evaluation of various lattice sums, and many related topics.
The paper~\cite{rpb269} discusses a recently-discovered connection
between the Gauss-Legendre algorithm and two of the Borwein algorithms
for computing $\pi$.

Perhaps Jon's best-known number-theoretic result, and one in whose
discovery experimental mathematics played an important role, is his paper
\cite{B5} (with his brother Peter) on a cubic counterpart of Jacobi's
identity, and a corresponding analogue of the arithmetic-geometric mean
of Legendre and Gauss. Briefly, if $AG_2(a,b)$ denotes the classical
arithmetic-geometric mean of $a$ and $b$, and $AG_3(a,b)$ denotes
the Borweins' cubic analogue, then
\[
AG_3(1,k) = {_2}F_1(1/3, 2/3; 1; 1-k^2)^{-1},
\]
which may be compared with Gauss's well-known identity
\[
AG_2(1,k) = {_2}F_1(1/2, 1/2; 1; 1-k^2)^{-1}.
\]
Here ${_2}F_1$ denotes a hypergeometric function.
The results given in~\cite{B5} 
have connections with the work of Ramanujan~\cite{Ramanujan}
on modular equations 
(see \cite{billion-digits} and the review of \cite{B5} by Bruce Berndt in 
\emph{Math.\ Reviews}).

\subsection{Mathematical Finance}			\label{subsec:finance}

Although Jon was primarily a pure mathematician, his interests extended much
further and included aspects of applied mathematics and statistics. This can
be illustrated by his work on mathematical finance. A
significant contribution is his 2014 paper~\cite{finance1} 
(with Bailey, de Prado and Zhu), provocatively titled
``Pseudo-mathematics and financial charlatanism''.

According to Bailey and Zhu \cite[pg.~229]{JBCC-proceedings}, 
the paper~\cite{finance1} grew
out of the authors' concern that, although mathematics had become a standard 
language to quantify financial phenomena, it was often used in a misguided
fashion, lending a patina of rigour to the topic at hand, but masking some
serious deficiencies. The paper demonstrated that many financial strategies
and fund designs, claimed to be backed by extensive ``backtests'' (analyses
based on historical market data), were nothing more than illusory artifacts
resulting from statistical overfitting. The conclusion was that backtest
overfitting is a likely reason why so many financial strategies and
fund designs fail, despite looking good on paper.
The paper~\cite{finance1} was followed by a more detailed analysis of how to
estimate the probability of backtest overfitting~\cite{finance2}.

Jon urged his co-authors to start a blog to present the conclusions of
\cite{finance1,finance2} to a broad audience. This resulted in the
\emph{Mathematical Investor} blog, whose mission is to identify and draw
attention to abuses of mathematics and statistics in the financial field,
and to call out the failure of many in the financial mathematics
community for silence on these abuses \cite[pg.~7]{JBCC-proceedings}.

\subsection{Mathematical Education and Visualisation}	\label{subsec:edu}

Jon had a passion for sharing his joy of mathematics
with students and a more general audience. Some of his research topics were
particularly well-suited for communication to high school and undergraduate
students. For example, we have already mentioned his interest in
experimental mathematics (\S\ref{subsec:expt}) and
$\pi$ (\S\ref{subsec:number_theory}). Jon was a keen blogger, starting this
activity in 2009 when he and David Bailey founded the ``Math Drudge'' 
(now ``Math Scholar'') blog,%
\footnote{\url{https://www.mathscholar.org}.}
which now has over $200$ articles on a wide range of topics, covering many
facets of modern mathematics, computing, and science.  For a sample of the
topics covered, see~\cite[xvii]{JBCC-proceedings}.

Jon regarded visualisation as a powerful tool for both experimental
mathematicians and mathematical communicators.  He was fond of the
maxims \emph{``it is often easier to see something than to explain it in
words''} and \emph{``a picture is worth a thousand words''}.  Examples may
be found in his paper~\cite{B10}, which contains striking images of walks in
the plane associated with various mathematical constants such as $\pi, e,$
and Champernowne's number $C_4$. In a different context, Jon used
visualisation to explain the dynamics of optimisation algorithms
such as Douglas-Rachford, for which see~\cite{Tam,Douglas-Rachford}.

\begin{figure}[ht]
\begin{center}
\includegraphics[width=16em]{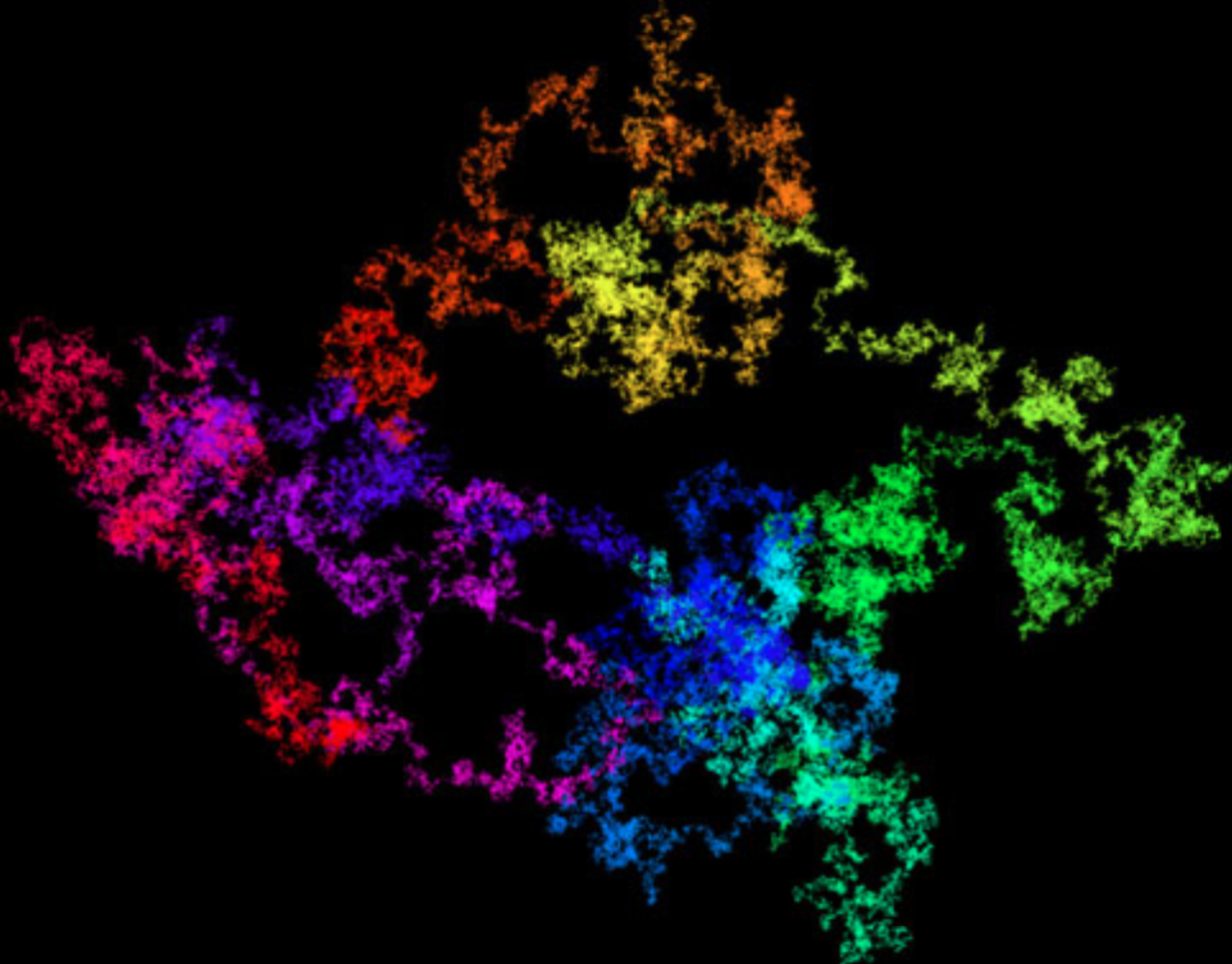}
\caption{A walk on the base-$4$ digits of $\pi$.} 
\end{center}
\end{figure}
\vspace*{-20pt}

A high-resolution version of this image, by Francisco
Arag\'on Artacho, is available at
\url{http://gigapan.com/gigapans/106803}.

\subsection*{Acknowledgements}

In writing this paper I have drawn heavily on information
provided in the JBCC Conference Proceedings by
Brailey Sims \emph{et al.}~\cite{JBCC-proceedings},
and the obituaries
by David Bailey~\cite{JBM-expt-math-blog}
and by Judith Borwein, Naomi Borwein, 
and Brailey Sims~\cite{obituary-AustMS}.
I have also made use of the extensive list of Jon's publications
maintained by David Bailey and Nelson Beebe~\cite{Beebe},
and information provided on the websites 
\cite{JB-homepage,JB-memorial,JBCC-website,JB-wikipedia}.
Thanks go to Heinz Bauschke, Tony Guttmann, Judy-anne Osborn, and Brailey Sims, 
for helpful comments on previous drafts.

\pagebreak[3]

\section{Annotated bibliography}			\label{sec:biblio}
		
An online list of Jon's publications is available
at~\cite{Beebe}. It contains over $2,000$ items
(although some are essentially duplicates), including (at least)
$387$ papers in refereed
journals, 
$103$ refereed or invited conference proceedings and book chapters,
and $26$ books.	
We list only a small sample, which includes Jon's better-known
books and most influential papers.

In the following, papers~\cite{B1}--\cite{Douglas-Rachford} 
include some of Jon's most
influential or highly cited papers. The list includes papers
that were given by Jon himself as
``ten career-best publications'' (up to mid-2013)
in an ARC grant application, and the comments on these papers
are  based on Jon's own.
Other papers have been selected by the author.

Jon wrote at least $26$ books. A few, like~\cite{B12},
are major research monographs. Some, like~\cite{Convex-optimization-book},
are advanced textbooks. Others are
``dictionaries'' (e.g.~\cite{B13,B14}) or edited collections
of papers by various authors on a particular topic
(e.g.~\cite{Pi-next-gen}).
Last but not least, 
there are some significant ``popular'' books written for a general
mathematical audience, e.g.\ those such as~\cite{crucible} designed to
popularise experimental mathematics.
We have included at least one from each category in the list
\cite{B13}--\cite{Pi-next-gen} below.

Finally, at the end, we give other relevant
sources~\cite{JBM-expt-math-blog}--\cite{JB-wikipedia}.

The papers and books by Jon are ordered chronologically
by date of publication; the final section is ordered alphabetically by author.

\small{

\vspace*{-10pt}
\renewcommand\refname{{\small Selected papers by Jonathan Borwein}}

 } 


\begin{thebibliography}{99}
\bibitem{rpb269}
Richard Brent,
The Borwein brothers, $\pi$ and the AGM,
in~\cite{JBCC-proceedings}, 323--348.

Discusses in detail some aspects of the book~\cite{B12} and proves that
the Gauss-Legendre algorithm for computing $\pi$
is equivalent to the Borwein-Borwein quartic algorithm.
\end{thebibliography}
\end{document}